# MAINTAINING DRAGGING AND THE PIVOT INVARIANT IN PROCESSES OF CONJECTURE GENERATION

Samuele Antonini

Department of Mathematics
University of Pavia (Italy)

Anna Baccaglini-Frank

Department of Mathematics
"Sapienza" University of Rome (Italy)

*In this paper, we analyze processes of conjecture generation in the context of open problems proposed in a dynamic geometry environment, when a particular dragging modality, maintaining dragging, is used. This involves dragging points while maintaining certain properties, controlling the movement of the figures. Our results suggest that the pragmatic need of physically controlling the simultaneous movements of the different parts of figures can foster the production of two chains of successive properties, hinged together by an invariant that we will call pivot invariant. Moreover, we show how the production of these chains is tied to the production of conjectures and to the processes of argumentation through which they are generated.*

## CONCEPTUAL FRAMEWORK AND RATIONALE

Dynamic Geometry Environments (DGEs) have acquired great interest for researchers in mathematics education and for teachers over the past years (e.g., Laborde & Strässer, 1990; Noss & Hoyles, 1996; Arzarello et al., 2002). In DGEs the figures can actually be seen and acted upon and their properties can be explored through dragging. This makes DGEs ideal for fostering, observing and analyzing processes of conjecture generation. In fact, the study of processes of conjecture generation, of argumentation and of proof in DGEs has come to be one of the leading themes of research in mathematics education (e.g., De Villiers, 1998; Hadas, Hershkowitz, & Schwarz, 2000; Arzarello et al., 2002; Mariotti, 2006).

The aim of the research we are presenting here is to yield a theoretical contribution for analyzing processes involved in conjecture generation and in argumentation within DGEs.

Dragging assumes a central role in the interaction with DGE figures, and various researchers have explored dragging modalities used by students during their explorations. In particular, Arzarello et al. (2002) and Olivero (2002) described different modalities used by the students according their goals, when generating conjectures in a DGE. Basing her work on such research, Baccaglini-Frank has identified and analysed students' use of four dragging modalities: free dragging, maintaining dragging, dragging with the trace mark active, dragging test (Baccaglini-Frank, 2010; Baccaglini-Frank & Mariotti, 2010). In this paper, we consider *maintaining dragging*, that consists in dragging some point intentionally maintaining invariant a certain property of the figure. This type of dragging is used especially in tasks that involve figuring out under which conditions certain properties are verified. In these cases, during the processes of conjecture generation a key role is played by the solvers' perception of invariants (Baccaglini-Frank, Mariotti & Antonini, 2009; Baccaglini-Frank, 2012; Leung, Baccaglini-Frank & Mariotti, 2013), described by Neisser (1989) as "aspects of stimulus information that persist despite movements". The perception of invariants involves the sensory experience of the solver in a DGE and it is distinct from the geometrical interpretation of the objects and of their mutual relations conceived within the theory of Euclidean geometry, a process that in the literature is referred to as *discernment* (Leung et al., 2013; Leung, 2008).

With the aim of clarifying cognitive processes in a DGE, Lopez Real and Leung (2006) distinguish between the realm of DGEs and the realm of Euclidean Geometry:

> A Dynamic Geometry Environment (DGE) is a computer microworld with Euclidean geometry as the embedded infrastructure. In this computational environment, a person can evoke geometrical figures and interact with them […]. It is a virtual mathematical reality where abstract concepts and constructions can be visually reified. In particular, the traditional deductive logic and linguistic-based representation of geometrical knowledge can be re-interpreted, or even redefined, in a DGE as dynamic and process-based interactive 'motion pictures' in real time. […] There appears to be a tension (rooted in the Euclidean view of what geometry is) pulling DGE research towards the direction of bridging an experimental–theoretical gap that seems to exist between the computational microworld and the formal abstract conceptual world.
>
> (Lopez Real & Leung, 2006, pp. 665-666)

In the realm of the theory of Euclidean geometry, geometrical properties are part of a network of logical relations, each validated by a mathematical proof.

On the other hand, in the realm of a DGE, points can be moved and the figures are modified as a consequence of such induced movement. The movement of a point can be direct (when the point itself is being dragged) or indirect (if the movement is a consequence of the dragging of another point): in this second case, the movement of the dragged point causes the movement of other points (because all the elements of the figure have to maintain the logical relationships imposed by the construction steps). The invariants are perceived in space but also in time and properties can be perceived simultaneously or in distinct temporal instances.

## CONTROLLING THE DGE FIGURE

To ease the reading of this paper we introduce the theoretical notions referring to one of the tasks assigned during the study.

Task: Construct: a point P and a line *r* through P, the perpendicular line to *r* through P, C on the perpendicular line, a point A symmetric to C with respect to P, a point D on the side of *r* containing A, the circle with center C and radius CP, point B as the second intersection between the circle and the line through P and D. Make conjectures about the possible types of quadrilateral ABCD can become, describing all the ways you can obtain a particular type of quadrilateral.

The figure resulting from the construction (Fig. 1) can be acted upon by dragging points (in our example we can think about dragging D), and some properties can be recognized as invariants for any movement of the dragged point (e.g., "CP = PA").

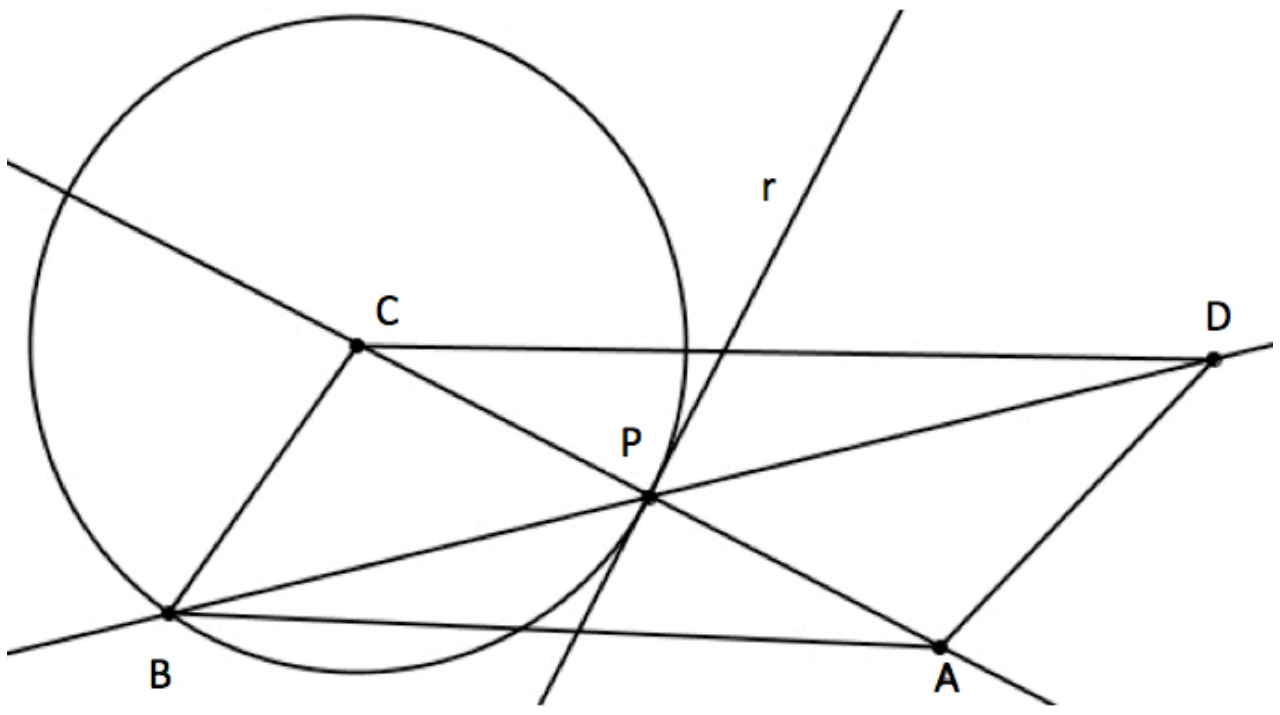

Figure 1: a possible result of the construction in the task above.

As the solver acts on the figure with the aim of generating a conjecture, s/he can decide to intentionally induce an invariant by dragging a point (e.g., "ABCD parallelogram" by dragging D), performing *maintaining dragging*.

The use of maintaining dragging involves dragging a point maintaining certain properties as invariants. Moving a point so that a DGE figure maintains a certain property requires a high level of control over the movement of different parts of the figure. The solver has to manage the relationships between the movements of the different parts of the figure (a similar case of coordinating movements has been studied in relation to cognitive processes involved in the use of pantographs, see Martignone & Antonini, 2009). The solver usually exercises indirect control over the invariant to maintain: its movement depends on the movement of the dragged point, that the solver controls directly, and sometimes the movements of the different parts of the figure are difficult to coordinate. The reader can experience the difficulty by trying to drag D (as in Fig. 1) maintaining the property "ABCD parallelogram".

To maintain an interesting configuration (A) the solver needs to control the figure more directly, so s/he passes from property A to a property $A_1$, which is easier and more direct to control, and such that its presence guarantees the presence of A. The nature of this process is abductive, similar to that described in Arzarello et al. (2002), were for *abduction* the following type of inference is intended: (*fact*) a fact A is observed; (*rule*) if C were true, then A would certainly be true; (*hypothesis*) so, it is reasonable to assume C is true (Pierce, 1960). In our case, given $A_1$, a property $A_2$ is generated such that $A_2 \Rightarrow A_1$. The spark initiating this abductive process is the need of better controlling the movements of the figure. Iterating the process may lead to an abductive chain of properties $A_1$, $A_2$, …$A_n$, such that each property ensures better control over the desired configuration A. This chain can later, in the proving phase, be flipped into a chain of deductive implications ($A_n \Rightarrow \ldots \Rightarrow A_2 \Rightarrow A_1$). So the spark initiating the successive inferences finds its origin in the pragmatic need of controlling the figure, and in particular in the need to coordinate the movement of the dragged point with that of the other parts of the figure in order to maintain invariant a desired property.

We stress how the search for a logical relation within the theory of Euclidean geometry sparks from an experience within the realm of the DGE. The solver will seek for fragments of theory with the goal of ameliorating his/her haptic control over the figure. These same fragments of theory can later be used for constructing an argumentation and finally a proof for the conjecture s/he will have reached.

**The pivot invariant**

During a first phase of the exploration leading to a conjecture, the chain $A_n \Rightarrow \ldots \Rightarrow A_2 \Rightarrow A_1 \Rightarrow A$ is developed as the solver searches for invariants to control more easily during maintaining dragging. When maintaining dragging is used, in a second phase, a new invariant $B_1$ can be perceived simultaneously, as $A_n$ is maintained. For the solver this new property $B_1$ has a very different status than the properties $A_i$. First of all, $B_1$ can be controlled directly. Secondly, the relation between the properties $B_1$ and $A_n$ is frequently perceived as causal (hence the arrow "$\rightarrow$" and not the implication symbol "$\Rightarrow$" in Figure 2) in the realm of DGE: the presence of $B_1$ guarantees (visually, for now) the simultaneous presence of $A_n$. Later, this causal relation can be interpreted logically as the implication $B_1 \Rightarrow A_n$. The process of conjecture generation may continue, leading to a second invariant $B_2$ simultaneously perceived during dragging, a third one $B_3$, etc., up to the generation of $B_m$ . These geometrical properties form a new chain of relationships perceived as causal relations in the realm of the DGE (Figure 2). The invariant $A_n$ acts as a pivot between the two chains of invariants and plays a fundamental role in the development of a conditional link between the properties that become the premise ($B_m$) and the conclusion (A) of the conjecture generated. We call $A_n$ the *pivot invariant*.

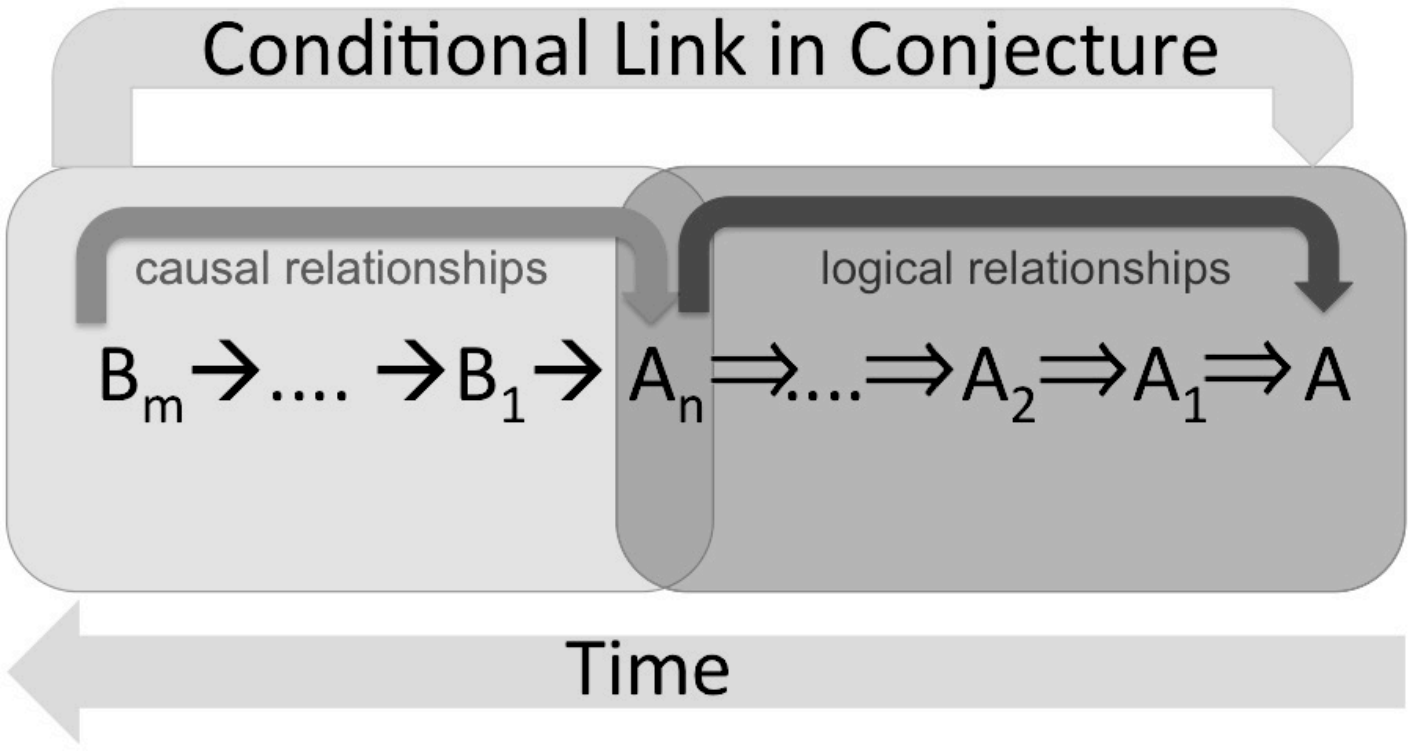

Figure 2: the two chains of invariants linked by the pivot invariant, discovered during the exploration that proceeds in time from right to left.

## THE PIVOT INVARIANT AT WORK: ANALYSIS OF A CASE

To show how the notion of pivot invariant can bring insight to students' conjecturing activity in a DGE we now present excerpts from an interview conducted with two students, Ste and Giu, both in the second year of high school (15 and 16 years old, respectively) in a northern Italian Liceo Scientifico. They had used Cabri II Plus the year before and they had been introduced to the four types of dragging introduced earlier in this paper (free dragging, maintaining dragging, dragging with the trace mark active, dragging test) during two previous lessons (for a complete description of the study see Baccaglini-Frank, 2010, or Baccaglini-Frank & Mariotti, 2010). The excerpts we present are from the part of the interview on the problem presented as an example in the previous section of this paper. The students' exploration lasted approximately 20 minutes, from the construction of the figure to the writing of their first conjecture on the possibility of obtaining a parallelogram. We analyze the excerpt in which the students identify the parallelogram as a possible configuration and decide to try to maintain it (the name of the student holding the mouse is in bold).

Giu: …So to get a parallelogram, what it looked like in the beginning…

**Ste**: So BP = PD by definition, that we know for sure.

Giu: Yes, yes, because they intersect at their midpoints, they are the diagonals.

Ste and Giu write a first conjecture: "ABCD is a parallelogram when BP = PD (that is when P is the midpoint of BD)". As Ste writes, Giu thinks aloud and proposes a proof for the conjecture.

Giu: Why? Because when this here is a diagonal [BD] of the parallel… of the quadrilateral, this here is another diagonal [CA] of the parall…of the quadrilateral.

Ste: But I need to add [referring to the written conjecture] in a parenthesis that CP is equal to PA.

Giu: CP = PA by definition, BP is equal to PD because we said so…and so they are diagonals that intersect at their midpoints, so it is a parallelogram.

The students have produced the following chain of deductions: BP = PD ($A_2$) ⇒ diagonals intersect at their midpoints ($A_1$) ⇒ ABCD is a parallelogram (A).

At this point the students start searching for ways of acting upon the figure to maintain property A during dragging. They seem to be seeking for ways to better control the movement of the different parts of the figure, that they still seem to find difficult to control through $A_2$. Suddenly Giu constructs the circle $C_{PD}$ (a circle with center in P and radius PD, see the Fig. 3) and says:

**Giu**: But see, you can do it like this. You can see that like this is comes out only when…no, you see…[he drags D so that $C_{PD}$ passes through the intersection, defined as B, between the line PD and the circle $C_{CP}$] […] You try and see *so that* this thing [concurrence of $C_{PD}$, $C_{CP}$ and *t*]…is maintained. (Figure 3)

Ste takes back the mouse and the students' attention shifts to the passage of $C_{PD}$ through B, the intersection between the line *t* (through P and D) and $C_{CP}$ (Figure 3).

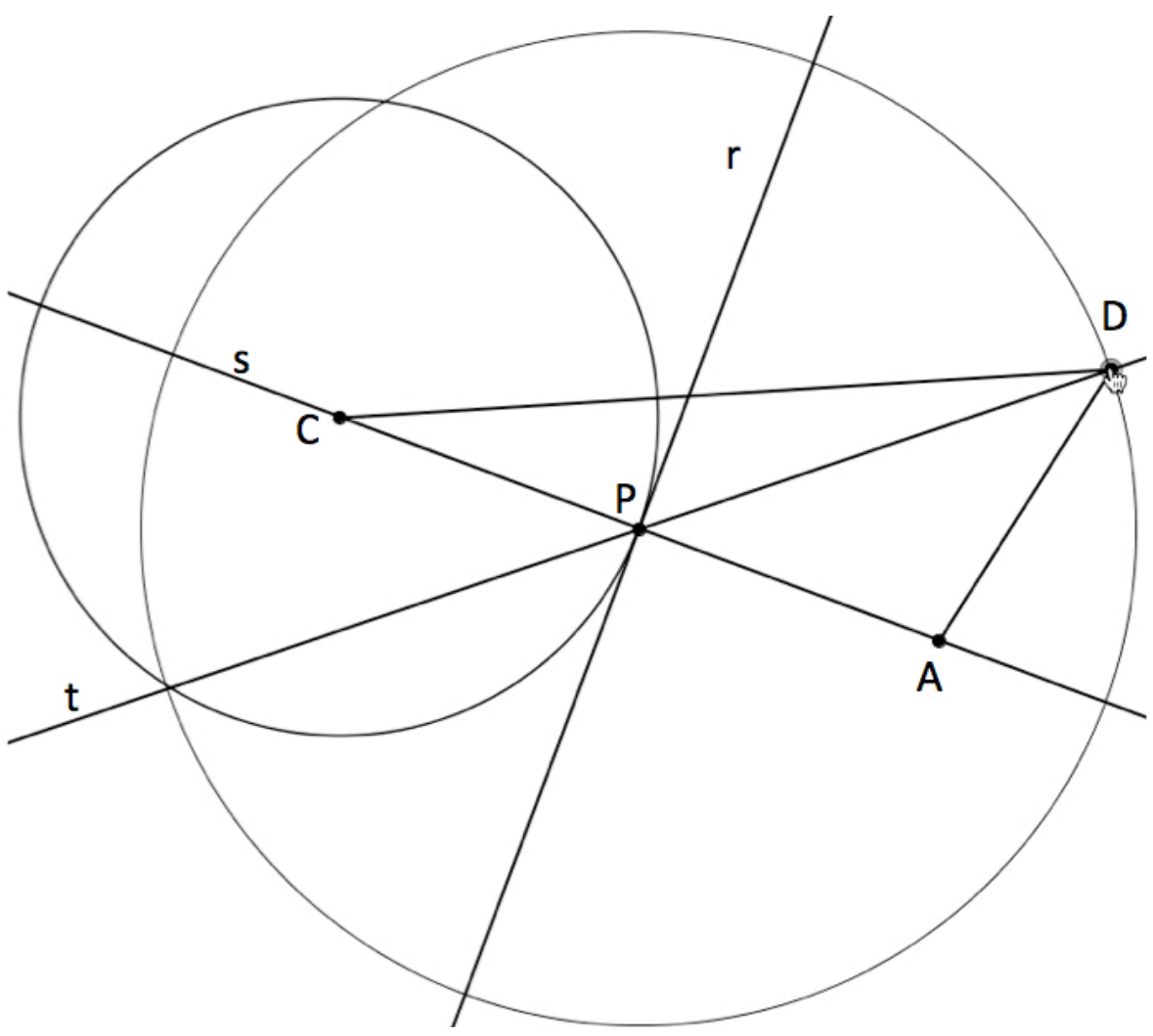

Figure 3: The students perform maintaining dragging inducing $C_{PD}$ to pass through the intersection of $C_{CP}$ with line *t*, defined as B in the original construction.

This is a new property $A_3$ that the students try to induce as an invariant. This property has been inferred through an abduction:

(*facts*) BP = PD ($A_2$), B is the second intersection of *t* with $C_{CP}$.

(*rule*) If B lies on $C_{PD}$, then BP = PD.

(hypothesis) B lies on $C_{PD}$ ($A_3$).

Dragging D it seems to be easier for the students to control $A_3$ than $A_2$. This is a fundamental moment in the process of conjecture generation: the invariant $A_3$ is discerned through an abduction in order to better control the figure and it is used for performing maintaining dragging in search for new invariants that might be causing $A_3$ (and the chain of invariants $A_3 \Rightarrow A_2 \Rightarrow A_1 \Rightarrow A$) to be visually verified. Elements of different natures concur in the discernment of this invariant: some are theoretical (e.g., the "rule" in the abduction is a theorem of Euclidean geometry) and others related to the phenomenology of the DGE (also see Leung et al., 2013).

As Giu continues to explore "when" $A_3$ is visually verified, he asks Ste to take back the mouse to concentrate on the movement of D when the trace is activated on it.

Giu: You maintain these things [B on $C_{PD}$]…it looks like a curve.

**Ste**: It's really hard!

Giu: Yes, I know…I can only imagine. It looks like a circle…with center in A.

**Ste**: It has to necessarily have radius AD! Anyway you would need AP to equal AD [he holds the mouse but stops dragging].

Ste is concentrated on maintaining the invariant $A_3$ while Giu tries to geometrically describe the trace mark. The students discern two invariants during dragging: D ∈ $C_{AP}$ ($B_1$), and PA = AD ($B_2$). Once the students construct $C_{AP}$, perform a dragging test dragging D along it, and notice that ABCD does seem to remain a parallelogram in this case, they write their final conjecture: "ABCD is a parallelogram if PA = AD" ($B_2 \Rightarrow A$).

The invariants $A_i$ are interpreted as logical consequences of the invariants $B_j$, and they are explicitly linked to them through the invariant $A_3$, which is a pivot-invariant. While the invariants $A_i$ arise mostly thanks to the students' theoretical knowledge of Euclidean geometry, the invariants $B_1, \ldots, B_m$ appear thanks to support offered within the phenomenological domain of the DGE, where invariants can be perceived through simultaneous perception accompanied by different levels of pragmatic control over the varying parts of the figure. The theoretical and pragmatic domain are hinged together by the pivot invariant that comes to life with hybrid characteristics.

## CONCLUSIONS

The analysis provided in this paper clarifies a specific process of conjecture generation in a DGE. In particular, the pragmatic need of physically controlling the figures explains how, in certain cases, the search for logical relations in the theory of Euclidean geometry can be fostered, together with the production of a chain of abductions leading to the conjecture. This new way of explaining how students come to substitute an invariant to maintain with a new property generalizes and ameliorates earlier descriptions of the process (see Baccaglini-Frank, 2010; Baccaglini-Frank & Mariotti, 2010): the process is now explained through the necessity of better controlling the figure, which leads to the production of two abductive chains.

The study points to at least two new directions of research. One is theoretical: the two chains of inferences may not necessarily be produced one after the other. In general, the two chains might be intertwined, leading to a greater complexity that needs to be further investigated. The second direction is practical and involves teachers. We believe that the theoretical notion of pivot invariant could be useful for a teacher who decides to promote students' conjecture generation in a DGE. Indeed, s/he could use it to gain deeper insight into students' processes of conjecture generation, and thus to better guide students' processes of conjecturing, argumentation, and proof.

### References

Arzarello, F., Olivero, F., Paola, D., & Robutti, O. (2002). A Cognitive Analysis of Dragging practises in Cabri Environments. *ZDM*, *34*(3), 66–72.

Baccaglini-Frank, A., Mariotti, M. A., & Antonini, S. (2009). Different Perceptions of Invariants and Generality of Proof in Dynamic Geometry. In Tzekaki, M., & Sakonidis, H. (Eds.), *Proc. of the 33rd Conference of the International Group for the Psychology of Mathematics Education,* Vol. 2, pp. 89-96. Thessaloniki, Greece: PME

Baccaglini-Frank, A. (2010). Conjecturing in Dynamic Geometry: A Model for Conjecture-generation through Maintaining Dragging. *Doctoral dissertation*, University of New Hampshire, Durham, NH.

Baccaglini-Frank, A. (2012). Dragging and Making Sense of Invariants in Dynamic Geometry. *Mathematics Teacher, 105 (8), 616-620.*

Baccaglini-Frank, A., & Mariotti, M.A. (2010) Generating Conjectures in Dynamic Geometry: the Maintaining Dragging Model. *International Journal of Computers for Mathematical Learning*, *15*(3), 225-253.

De Villiers, M. (1998). An alternative approach to proof in dynamic geometry. In R. Leher & D. Chazan (Eds.), *Designing Learning Environments for Developing Understanding of Geometry and Space* (pp. 369-393). Hillsdale, NJ: Lawrence Erlbaum Associates.

Hadas, N., Hershkowitz, R., & Schwarz, B. B. (2000). The Role of Contradiction and Uncertainty in Promoting the Need to Prove in Dynamic Geometry Environments. *Educational Studies in Mathematics, 44*(1/2), 127-150.

Laborde, J. M., & Strässer, R. (1990). Cabri-Géomètre: A microworld of geometry for guided discovery learning. *ZDM, 22*(5), 171–177.

Leung, A. (2008). Dragging in a Dynamic Geometry Environment Through the Lens of Variation. *Intern. Journal of Computers for Mathematical Learning*, *13*(2), 135–157.

Leung, A., Baccaglini-Frank, A. & Mariotti, M.A. (2013). Discernment in dynamic geometry environments. *Educational Studies in Mathematics*, 84(3), 439–460.

Lopez-Real, F. & Leung, A. (2006). Dragging as a conceptual tool in dynamic geometry environments. *International Journal of Mathematical Education in Science and Technology, 37* (6), 665-679.

Mariotti M.A. (2006) Proof and proving in mathematics education. A. Gutiérrez, & Boero, P. (Eds.) *Handbook of Research on the Psychology of Mathematics Education* (pp. 173-204). Sense Publishers, Rotterdam, The Netherlands.

Martignone, F. & Antonini, S. (2009). Exploring the mathematical machines for geometrical transformations: a cognitive analysis. In *Proc. of the 33rd Conference of the International Group for the Psychology of Mathematics Education,* Vol.4, pp. 105-112. Thessaloniki, Greece: PME.

Neisser, U. (Ed.). (1989). *Concepts and conceptual development: Ecological and intellectual factors in categorization*. Cambridge, UK: Cambridge University Press.

Noss, R. & Hoyles, C. (1996) *Windows on Mathematical Meanings Learning Cultures and Computers*, Dordrecht, The Netherlands: Kluwer Academic Publishers.

Olivero, F. (2002). *The Proving Process within a Dynamic Geometry Environment*. PhD Thesis, University of Bristol.

Peirce C. S. (1960). *Collected Papers II, Elements of Logic*. Harvard: University Press.